\documentclass[12pt]{amsart}
\usepackage{amsfonts,amsmath,amsthm,amssymb}
\usepackage{latexsym}
\usepackage{enumerate}
\usepackage{graphics}
\newtheorem{theorem}{Theorem}[section]

\theoremstyle{definition}

\theoremstyle{remark}

\theoremstyle{remark}

\newcommand{\beql}[1]{\begin{equation}\label{#1}}
\newcommand{\eeq}{\end{equation}}

\begin{document}

\title{A note on the existence of $BH(19,6)$ matrices}

\author{Ferenc Sz\"oll\H{o}si}

\date{December, 2011.\\ Preprint}

\address{Ferenc Sz\"oll\H{o}si: Department of Mathematics and its Applications, Central European University, H-1051, N\'ador u. 9, Budapest, Hungary}\email{szoferi@gmail.com}


\dedicatory{Dedicated to Professor Kathy J. Horadam on the occasion of her $60$th birthday}

\begin{abstract}
In this note we utilize a non-trivial block approach due to Petrescu \cite{MP1} to exhibit a Butson-type complex Hadamard matrix of order $19$, composed of sixth roots of unity.
\end{abstract}

\maketitle

	

\section{A new Butson-type complex Hadamard matrix}

A complex Hadamard matrix $H$ is an $n\times n$ complex matrix with unimodular entries such that $HH^\ast=nI_n$, where $\ast$ denotes the Hermitian adjoint and $I_n$ is the identity matrix of order $n$. Throughout this note we are concerned with the special case when the entries are some $q$th roots of unity. These matrices are called Butson-type complex Hadamard matrices and are denoted by $BH(n,q)$ \cite{AB1}. The existence of $BH(n,q)$ matrices is wide open in general, but it is believed that real, i.e. $BH(4k,2)$ matrices exist for every integral number $k$. This is the famous Hadamard conjecture. Butson-type complex Hadamard matrices have applications in signal processing, coding theory \cite{HLT} and harmonic analysis \cite{KM2}, among other things. They can also lead to constructions of real Hadamard matrices. This approach was demonstrated very recently in \cite{CCL}, where $BH(n,6)$ matrices were considered. For further results on $BH(n,q)$ matrices we refer the reader to Horadam's celebrated book \cite{KJH1}.

The main result of this note is the following
\begin{theorem}\label{thm1}
There exists a $BH(19,6)$ matrix.
\end{theorem}
Order $n=19$ was listed as the smallest outstanding order of $BH(n,6)$ matrices in \cite{CRC1}. Addressing the existence of complex Hadamard matrices of prime orders is a notoriously difficult problem. One of the reasons for this is that the standard construction methods developed for the study of real Hadamard matrices and various orthogonal arrays fundamentally rely on ``plug-in'' methods and block constructions resulting in objects of composite order.

However, a non-standard block approach was utilized by Petrescu in his PhD thesis \cite{MP1}, who exhibited various complex Hadamard matrices of prime orders. We recall his method as follows. Let $s$ be a positive integer, $X, Y$ be $s\times s$, $D$ be an $(s+1)\times (s+1)$ matrix, and assume that the $s\times (s+1)$ matrix $T$ is composed of $s$ noninitial rows of a normalized complex Hadamard matrix of order $s+1$. The following matrix of order $3s+1$
\[H=\left[\begin{array}{ccc}
X & Y & T\\
Y & X & T\\
T^\ast & T^\ast & D\\
\end{array}
\right]\]
is called Petrescu's array. The matrix $H$ is complex Hadamard if and only if the entries of $X,Y,D$ and $T$ are unimodular and the following orthogonality equations are satisfied:
\begin{eqnarray*}
2T^\ast T+DD^\ast&=&(3s+1)I_{s+1},\label{1}\\
XX^\ast+YY^\ast+TT^\ast&=&(3s+1)I_s,\label{2}\\
(X+Y)T+TD^\ast&=&0,\label{3}\\
XY^\ast+YX^\ast+TT^\ast&=&0.\label{4}
\end{eqnarray*}
Let us denote by $J_{s+1}$ the all-$1$ matrix of order $s+1$. The system of equations above can be transformed into the following equivalent, but more informative one (see \cite{MP1} for details):
\begin{eqnarray}
DD^\ast=D^\ast D&=&(s-1)I_{s+1}+2J_{s+1},\label{5}\\
DJ&=&JD,\label{6}\\
X+Y&=&-\frac{1}{s+1}TD^\ast T^\ast,\label{7}\\
(X-Y)(X-Y)^\ast&=&(3s+1)I_s.\label{8}
\end{eqnarray}
Thus the lower right corner $D$ is a normal, regular matrix, and hence in view of \eqref{5} can be thought as some kind of generalized, ``complex biplane''. The sum $X+Y$ is essentially a unitary transform of $D^\ast$, while the difference $X-Y$ is an orthogonal matrix.

In order to prove Theorem \ref{thm1} we take $s=6$ and look for some solutions of the above equations in which all of the matrix entries are some $6$th roots of unity. First we search for candidate matrices $D$ satisfying \eqref{5} and \eqref{6}, and then we search for a suitable, normalized $BH(7,6)$ matrix, whose six noninitial rows will constitute $T$ such that the right hand side of \eqref{7} is a matrix in which all entries are sums of some sixth roots of unity. Finally, we decompose the right hand side of \eqref{7} into some matrices $X$ and $Y$ and check whether they meet the final condition \eqref{8}. After implementing this rather straightforward algorithm in {\it Mathematica}, we have found the following matrix, virtually in seconds:
\[W_{19}=\left[
\begin{array}{cccccc|cccccc|ccccccc}
 3 & 0 & 1 & 1 & 0 & 0 & 5 & 4 & 3 & 5 & 3 & 2 & 1 & 1 & 3 & 5 & 4 & 3 & 0 \\
 0 & 0 & 1 & 3 & 3 & 1 & 4 & 2 & 4 & 5 & 1 & 5 & 1 & 4 & 3 & 3 & 1 & 5 & 0 \\
 0 & 0 & 1 & 4 & 2 & 4 & 2 & 4 & 3 & 2 & 4 & 1 & 3 & 3 & 1 & 4 & 5 & 1 & 0 \\
 1 & 2 & 4 & 2 & 1 & 2 & 4 & 4 & 2 & 4 & 5 & 0 & 3 & 5 & 1 & 1 & 3 & 4 & 0 \\
 2 & 5 & 4 & 3 & 2 & 0 & 4 & 2 & 0 & 1 & 4 & 2 & 4 & 1 & 5 & 3 & 1 & 3 & 0 \\
 0 & 3 & 5 & 4 & 5 & 0 & 4 & 5 & 3 & 1 & 3 & 4 & 5 & 3 & 4 & 1 & 3 & 1 & 0 \\
 \hline
 5 & 4 & 3 & 5 & 3 & 2 & 3 & 0 & 1 & 1 & 0 & 0 & 1 & 1 & 3 & 5 & 4 & 3 & 0 \\
 4 & 2 & 4 & 5 & 1 & 5 & 0 & 0 & 1 & 3 & 3 & 1 & 1 & 4 & 3 & 3 & 1 & 5 & 0 \\
 2 & 4 & 3 & 2 & 4 & 1 & 0 & 0 & 1 & 4 & 2 & 4 & 3 & 3 & 1 & 4 & 5 & 1 & 0 \\
 4 & 4 & 2 & 4 & 5 & 0 & 1 & 2 & 4 & 2 & 1 & 2 & 3 & 5 & 1 & 1 & 3 & 4 & 0 \\
 4 & 2 & 0 & 1 & 4 & 2 & 2 & 5 & 4 & 3 & 2 & 0 & 4 & 1 & 5 & 3 & 1 & 3 & 0 \\
 4 & 5 & 3 & 1 & 3 & 4 & 0 & 3 & 5 & 4 & 5 & 0 & 5 & 3 & 4 & 1 & 3 & 1 & 0 \\
 \hline
 5 & 5 & 3 & 3 & 2 & 1 & 5 & 5 & 3 & 3 & 2 & 1 & 0 & 0 & 0 & 0 & 1 & 1 & 3 \\
 5 & 2 & 3 & 1 & 5 & 3 & 5 & 2 & 3 & 1 & 5 & 3 & 0 & 0 & 1 & 3 & 0 & 1 & 0 \\
 3 & 3 & 5 & 5 & 1 & 2 & 3 & 3 & 5 & 5 & 1 & 2 & 1 & 3 & 0 & 0 & 0 & 1 & 0 \\
 1 & 3 & 2 & 5 & 3 & 5 & 1 & 3 & 2 & 5 & 3 & 5 & 0 & 1 & 0 & 1 & 0 & 4 & 1 \\
 2 & 5 & 1 & 3 & 5 & 3 & 2 & 5 & 1 & 3 & 5 & 3 & 1 & 0 & 4 & 1 & 1 & 0 & 0 \\
 3 & 1 & 5 & 2 & 3 & 5 & 3 & 1 & 5 & 2 & 3 & 5 & 1 & 0 & 1 & 0 & 4 & 0 & 1 \\
 0 & 0 & 0 & 0 & 0 & 0 & 0 & 0 & 0 & 0 & 0 & 0 & 4 & 1 & 1 & 0 & 1 & 0 & 0 \\
\end{array}
\right].\]
In the matrix $W_{19}$ above we have replaced the sixth roots of unity $\mathbf{e}^{\frac{2\pi\mathbf{i}}{6}j}$ with $j$, $j=0,1,\hdots,5$ for typographical reasons.
A more in-depth analysis of this approach along with a handful of additional examples of Petrescu-type complex Hadamard matrices will be presented elsewhere.
\section*{Acknowledgements}
The results of this note were obtained during a research visit to Professor Robert Craigen. The author acknowledges the warm hospitality of the Department of Mathematics at University of Manitoba. This work was supported by the Hungarian National Research Fund OTKA $K$-$77748$ and by the Doctoral Research Support Grant of CEU.

\end{document}